\documentclass[letterpaper, 10 pt, conference]{ieeeconf}  

\IEEEoverridecommandlockouts                              
\usepackage{times} 
\usepackage{amsmath,amsfonts} 
\usepackage{amssymb}  
\usepackage{tabularx}
\usepackage[normalem]{ulem}
\usepackage{dsfont}
\usepackage{graphicx}
\usepackage{subcaption}
\usepackage[dvipsnames,usenames]{xcolor}
\usepackage[colorlinks,%
linkcolor=BrickRed,%
filecolor=BrickRed,%
citecolor=RoyalPurple,%
]{hyperref}
\usepackage{algorithm}
\usepackage{algorithmic}
\usepackage{float}

\newcommand{\rebuttal}[1]{{ #1}}

\def\R{\mathbb{R}}

\def\notes#1{\marginpar{\tiny #1}\typeout{Notes!
Notes!
Notes!
}}
\renewcommand{\notes}[1]{\typeout{notes!}}

\def\R{\mathbb{R}}

\newtheorem{remark}{Remark}

\def\beq{\begin{eqnarray}} 
\def\bc{\begin{center}} 
\def\be{\begin{enumerate}}
\def\bi{\begin{itemize}} 
\def\bs{\begin{small}}
\def\bS{\begin{slide}}
\def\ec{\end{center}} 
\def\ee{\end{enumerate}}
\def\ei{\end{itemize}}
\def\es{\end{small}}
\def\eS{\end{slide}}
\def\eeq{\end{eqnarray}}

\newcounter{rmnum}

\newcounter{anum}

\title{\LARGE \bf
Fast filtering of non-Gaussian models using 
\\Amortized Optimal Transport Maps
}

\author{Mohammad Al-Jarrah$^{\star,\dagger}$, Bamdad Hosseini$^\dagger$, Amirhossein Taghvaei$^\star$
\thanks{Mohammad Al-Jarrah and Amirhossein Taghvaei are supported by the National Science Foundation (NSF) award EPCN-2318977. Bamdad Hosseini is supported by the NSF award DMS-2208535}
    {\thanks{$^\star$Department of Aeronautics \& Astronautics, University of Washington, Seattle; {\tt\small mohd9485@uw.edu,amirtag@uw.edu}.}}
    {\thanks{$^\dagger$Department of Applied Mathematics, University of Washington, Seattle
        {\tt\small mohd9485@uw.edu,bamdadh@uw.edu}.}}
}

\begin{document}

      \maketitle
	 \thispagestyle{empty}
	 \pagestyle{empty}


	 \begin{abstract}
In this paper, we present the amortized optimal transport filter (A-OTF) designed to mitigate the computational burden associated with the real-time training of optimal transport filters (OTFs). OTFs can perform accurate non-Gaussian Bayesian updates in the filtering procedure, but they require training at every time step, which makes them expensive. The proposed A-OTF framework exploits the similarity between OTF maps during an initial/offline training stage in order to reduce the cost of inference during online calculations.  More precisely, we use clustering algorithms to select relevant subsets of pre-trained maps whose weighted average is used to compute the A-OTF model akin to a mixture of experts.  A series of numerical experiments validate that A-OTF achieves substantial computational savings during online inference while preserving the inherent flexibility and accuracy of OTF.
	\end{abstract}

    \section{Introduction}
The nonlinear filtering problem concerns the approximation of the conditional probability distribution (posterior) of the hidden state of a stochastic dynamical system given a sequence of partial and noisy observations. This problem is crucial in a wide range of applications including: 
satellite orbit determination and  navigation systems~\cite{jazwinski2007stochastic}; 
weather forecasting~\cite{van2010geosciences}; 
machine learning~\cite{bishop2006pattern};  and
economics~\cite{brigo1998some,javaheri2003filtering}.


Traditional methods for nonlinear filtering, such as the Kalman filter (KF) and its extensions~\cite{kalman1960new,kalman1961new}, simplify the problem by imposing a Gaussian ansatz on the joint 
distribution of the state and observations. This  Gaussian assumption 
is adequate for nearly linear systems with additive Gaussian noise but it is ineffective for highly nonlinear and complex systems. On the other hand, particle filter methods such as sequential importance re-sampling (SIR)~\cite{gordon1993novel,doucet09} approximate the posterior 
with the weighted empirical distribution of a set of particles. 
As a result, particle filters are able to handle nonlinear and non-Gaussian models but they often suffer from the weight degeneracy issue, usually in high-dimensional settings \cite{bickel2008sharp,snyder2008obstacles}.

To overcome the aforementioned limitations, 
coupling and transport-based approaches have been developed ~\cite{reich2013nonparametric,marzouk2016introduction,spantini2022coupling,shi2022conditional,taghvaei2023survey}, with the optimal transport filter (OTF)~\cite{taghvaei2022optimal,al2023optimal,grange2023computational,al-jarrah2024nonlinear} emerging as a particularly attractive alternative. The OTF deterministically transports particles using  an optimal transport (OT) formulation of the prior to posterior relationship. OTF offers two main advantages: (i) it requires only samples from the joint distribution of states and observations, eliminating the need for explicit analytical models of the likelihood and dynamics, and (ii) it enables flexible neural network parameterizations to enhance the expressivity of transport maps. However, the computational burden of OTF is significant as the OT maps are 
computed online at each time step via a stochastic 
optimization procedure.


The goal of this article is to reduce the computational cost 
of OTF during online/real-time filtering by exploiting the 
shared structure of OTF maps across various time instances: simply put, similar state and observation distributions 
should lead to similar OTF maps and this ``similarity" 
can be exploited to improve efficiency of the algorithm.
Inspired by the amortized optimization methodology~\cite{amos2023tutorial,gershman2014amortized,rezende2014stochastic,kingma2013auto,stuhlmuller2013learning}, we propose the amortized OTF (A-OTF). 
This method consists of an offline and an online stage:
During the offline stage
a collection of OTF maps are computed from a set of 
(often random) simulations of the model.
During the online stage, the prior-to-posterior map is 
approximated as a weighted linear combination of the offline OTF 
maps. The defining weights for the A-OTF maps are determined 
using similarity measures between the current prior and those associated 
with the pre-trained maps, an approach analogous to kernel
interpolation or scattered data approximation~\cite{nadaraya1964estimating, wendland2004scattered}. To further accelerate the online computations and avoid redundancy in the pre-trained set, we also apply a clustering algorithm in the offline stage that 
allows us to localize the A-OTF map around the current prior. Specifically, we use the partitioning around medoids algorithm~\cite{kaufmann1987clustering,kaufman1990partitioning,kaufman2009finding}
which we refer to as the $K$-Medoids. As it is often the case with local interpolation/regression, 
the quality of A-OTF is tied to the quality of the pre-training data set and the 
complexity of the underlying dynamics. Indeed, in the previous work \cite{al2024data}, 
the authors utilized a single pre-trained OTF filter in real-time for stationary 
dynamics. The A-OTF can be viewed as an extension of this idea to  non-stationary
settings.

The rest of the paper is organized as follows: 
Sec.~\ref{sec:problem_formulation} includes the mathematical setup and the modeling assumptions; Sec.~\ref{sec:amortized} contains the proposed methodology; and section~\ref{sec:numerics} presents several numerical experiments 
and benchmarks.

    \section{Problem Formulation}\label{sec:problem_formulation}
    \subsection{Setup}
\label{sec:filtering-problem}
Consider the discrete-time stochastic dynamical system 
\begin{subequations}\label{eq:model}
\begin{align}\label{eq:model-dyn}
    X_{t} &\sim a_t(\cdot|X_{t-1}) , \quad X_0 \sim \pi_0 \\\label{eq:model-obs}
    Y_t &\sim h_t(\cdot|X_t)
\end{align}
\end{subequations}
for $t=1,2,\ldots $ where $X_t\in \mathbb{R}^n$ is the hidden state of the system, $Y_t \in \mathbb{R}^{n_y}$ is the observation, $\pi_0$ is the initial probability distribution of the state, $a_t( x'|x)$ is the transition kernel from $X_{t-1}=x$ to $X_t=x'$,  and $h_t(y|x)$ is the  likelihood of 
observing $Y_t=y $ given $X_t=x$. 


The filtering problem is to infer the conditional distribution of the state $X_t$ given the history of the observations $\{Y_1,\dots,Y_{t}\}$, that is, the distribution
\[
\pi_t := \mathbb{P}(X_t\in \cdot |Y_1,\dots,Y_{t} )
\]
often referred to as the {\it posterior} distribution. 


The posterior distribution $\pi_t$ can be expressed via a recursive equation which is fundamental to the design of filtering algorithms. In particular, by leveraging the dynamic and observation models, we introduce the following propagation and conditioning operators:
\begin{subequations}
    \begin{align}
    \label{eq:propagation}  
    \text{(propagation)}~ \pi \mapsto \mathcal A_t \pi &:= \int_{\mathbb R^n} a_t(\cdot|x) \pi(x)d x,\\
  \text{(conditioning)}~ \pi \mapsto \mathcal B_{t,y}(\pi) &:= \frac{h_t(y|\cdot)\pi(\cdot)}{\int_{\mathbb R^n} h_t(y|x) \pi(x)d x},
  \label{eq:Bayesian}   
\end{align}
\end{subequations}
for an arbitrary probability distribution $\pi$. The propagation operator $\mathcal A_t$ represents the update for the distribution of the state according to the dynamic model~\eqref{eq:model-dyn} and the 
conditioning operator $\mathcal B_{t,y}$ represents the Bayes' rule that carries out the conditioning according to the observation model~\eqref{eq:model-obs}.
Combining these operators yields an update for
for the posteriors~\cite{cappe2009inference}:
\begin{equation}\label{eq:bayes}
 \pi_{t} \;=\; \mathcal{B}_{t,Y_t}(\mathcal{A}_t \pi_{t-1}).
\end{equation}

\subsection{The Optimal Transport Filter}

The OTF algorithm is introduced in the recent work of the authors~\cite{al-jarrah2024nonlinear} as an algorithm for approximating the Bayesian update~\eqref{eq:Bayesian} with an OT map which transports particles from the prior to the posterior. Let us briefly recall how OTF works. 
For any joint distribution $P_{X,Y}$ consider the optimization problem 
\begin{align}\label{eq:new_loss}
    \max_{ f \in c\text{-Concave}_x}\,
    \min_{T \in \mathcal{M}(P_X \otimes P_Y)}\, J(f,T;P_{X,Y}),
\end{align}
where $P_X \otimes P_Y$ denotes the independence coupling 
of the $X$ and $Y$ marginals $P_X$ and $P_Y$,
$\mathcal{M}(P_X \otimes P_Y)$ is the set of  maps $\R^n \times \R^{n_y} \mapsto \R^n$ that are  $P_X \otimes P_Y$-measurable, and the set $c\text{-Concave}_x$ denotes the set of functions $f(x,y)$ on $\R^n \times \R^{n_y} \mapsto \R$ that are $c$-concave in their first variable $x$ everywhere\footnote{$f$ is $c$-concave iff $\frac{1}{2}\|\cdot\|^2-f$ is convex.}. The explicit form of the objective $J$ appears in~\cite[eq(8c)]{al-jarrah2024nonlinear}. Assuming $P_X$ is absolutely continuous with respect to  the Lebesgue measure, the optimization~\eqref{eq:new_loss} has a unique solution pair $(\overline{f},\overline{T})$ where 
\begin{equation*}
     \overline{T}(\cdot,y)_{\#}P_X(\cdot) = P_{X|Y}(\cdot|y) \quad \forall y,
\end{equation*}
where ``\(\#\)'' denotes the push-forward operator~\cite[Prop. 2.3.]{al-jarrah2024nonlinear}.
The map $\overline{T}$, produced by the above optimization procedure, is used to replace the Bayesian step in~\eqref{eq:bayes}. Specifically, if we let $P_{X,Y}(x,y) = (\mathcal{A}_t\,\pi_{t-1}(x)) h_t(y|x)$ and write $T_t$ for the solution to~\eqref{eq:new_loss} we have the identity
\begin{equation}\label{eq:map_defn}
    T_t(\cdot,\,y)_{\#}(\mathcal{A}_t\,\pi_{t-1}) \;=\; \mathcal{B}_{t,y}(\mathcal{A}_t\,\pi_{t-1}),\quad \forall y,
\end{equation}
which yields the equivalent posterior update rule
\begin{equation}\label{eq:map_T_true}
    \pi_{t} 
    \;=\; 
    T_t(\cdot,\,Y_t)_{\#}\,\mathcal{A}_t\,\pi_{t-1}.
\end{equation}


This new update rule can be numerically implemented using an ensemble of $N$ particles $\{X_{t|t}^i\}_{i=1}^N$. In this case, the update can be done in two 
steps:
\begin{subequations}
    \begin{align}  
    \label{eq:X_particle_law}
    &\text{(propagation)}  &&X^i_{t|t-1} \;\sim\; a_t(\cdot|X^i_{t-1|t-1}),\\
  &\text{(conditioning)} &&X^i_{t|t} \;=\; \widehat{T}_t(X^i_{t|t-1},Y_t), 
  \label{eq:posterior_particle_law}
\end{align}
\end{subequations}
where the map $\widehat{T_t}$ is obtained by solving the optimization problem~\eqref{eq:new_loss} where the joint distribution $P_{X,Y}$ is approximated by the empirical distribution of samples  $(X_{t|t-1}^i,Y_{t|t-1}^i)_{i=1}^N$ with $Y_{t|t-1}^i\sim h_t(\cdot|X_{t|t-1}^i) $ and  $Y_{1:t}=\{Y_1,\dots,Y_t\}$ are the true observations.  Furthermore, the function $f$ and the map $T$ are parameterized as neural networks as described in~\cite{al-jarrah2024nonlinear}.
Finally, we note that while OTF gives an OT characterization of 
the conditioning maps, these can be computed in various 
other ways, see for example ~\cite{spantini2022coupling,wang2023efficient,zeng2024ensemble}.

\subsection{Objective}


Let us now 
suppose that an arbitrary set of pre-trained OTF maps $\{\widetilde T_m\}_{m=1}^M$ are given and that for each map we also have the training data $(\widetilde X_{m}^i,\widetilde Y_m^i)_{i=1}^N$. The goal is to approximate the solution to~\eqref{eq:new_loss} for a new set of samples $(X^i,Y^i)_{i=1}^N$. Thus, 
we consider the problem:
\begin{align*}
    \text{Given:}\quad & \left\{\left(\widetilde X^i_m,\widetilde Y_m^i\right)_{i=1}^N,\widetilde T_m \right\}_{m=1}^M,    \\
    \text{Approx.:}\quad   & \textit{Solution $T$ to~\eqref{eq:new_loss},} \\
    &\textit{for new set of samples} \; \left(X^i,Y^i\right)_{i=1}^N. 
\end{align*}
Henceforth we refer to the given data above as the {\it the training 
data} for A-OTF. 


 \section{The Amortized optimal transport filter}\label{sec:amortized}

Towards solving the problem outlined above we propose to 
approximate the map $T$ as a weighted combination of the 
pre-trained maps $\widetilde{T}_m$ with  the weights  determined by a similarity measure between the current empirical  prior (given by the 
samples $\{ X^i\}_{i=1}^N$) and those associated with 
the pre-trained maps (given by the samples $\{\widetilde{X}^i_m\}_{i=1}^N$). 
Simply put, we approximate the mapping $\{ X^i\}_{i=1}^N \mapsto T$ 
by regressing it over the A-OTF training data.
To accelerate this procedure, avoid redundancy, and improve robustness, 
we apply a clustering algorithm to the training data which allows us to build ``local" approximations 
to $T$. We present the details of this procedure 
below. 



\subsection{Details of the Offline Stage}
We employ  the $K$-medoids algorithm \cite{kaufmann1987clustering} 
to split the training data 
into $K$ clusters. The $K$-medoids algorithm identifies $K$ medoids by minimizing the total distance between each data point and the medoid of its assigned cluster. Unlike centroids, which do not necessarily correspond to actual data points, medoids are selected from within the dataset. This procedure necessitates the computation of pairwise distances, which are encapsulated in a distance matrix $D\in \mathbb{R}^{M\times M}$.
Let us introduce the point cloud
\[
\widetilde S_m := \Bigl((\widetilde X^i_m,\widetilde Y^i_m)_{i=1}^{N}, \widetilde T_m\Bigr),\quad m=1,\dots,M,
\]
in the product space of distributions on $\R^n \times \R^{n_y}$ 
and transport maps. 
For any pair $\widetilde S_u$ and $\widetilde S_v$ we define the $u,v$ element of the distance matrix $D$ as
\begin{align}\label{eq:Dmn}
   D_{u,v} := d\Bigl(\widetilde S_u, \widetilde S_v\Bigr),
\end{align}
for an appropriate distance function or metric $d$. In this 
work we consider three choices for $d$ as outlined in Table~\ref{tab:d_notation}.
Applying $K$-medoids to $D$ yields a clustering of the training data into 
$K$ subsets with representative medoids  $\{S_k^\ast\}_{k=1}^K$.

\begin{table}[t]
    \centering
    \renewcommand{\arraystretch}{1.5}
    \begin{tabular}{|p{1cm}|p{6.5cm}|}
    \hline
        \textbf{$d_{W_2}$} & Wasserstein-2 ($W_2$) distance  
        between $X$-samples $\{X_u^i\}_{i=1}^N, \{X_v^i\}_{i=1}^N$  \\
    \hline
    \textbf{$d_{MMD}$} & Maximum mean discrepancy (MMD) between 
    $X$-samples $\{X_u^i\}_{i=1}^N, \{X_v^i\}_{i=1}^N$ \\
    \hline
    \textbf{$d_{T}$} & Averaged distance between transported particles:\newline
    \(\frac{1}{2N}\sum_{i=1}^N \Bigl(\,
          \bigl\|T_u(X_u^i,Y_u^{\sigma_i}) 
            - T_v(X_u^i,Y_u^{\sigma_i})\bigr\|
           +\;\bigl\|T_u(X_v^i,Y_v^{\sigma_i}) 
            - T_v(X_v^i,Y_v^{\sigma_i})\bigr\| \Bigr)\) \newline
            where $\{\sigma_1,\sigma_2,\ldots,\sigma_N\}$ denotes an independent random permutation of the index set $\{1,2,\ldots,N\}$, and $\|\cdot\| $ is the Euclidean norm. 
    \\
    \hline
        \end{tabular}
    \caption{Choices of the distance function $d$ in $K$-medoids.}
    \label{tab:d_notation}
\end{table}


\subsection{Details of the Online Stage}
 At this stage, we approximate the map $T_t$ in~\eqref{eq:map_T_true} for a new ensemble of particles $(X^i_{t|t-1},Y^i_{t|t-1})_{i=1}^{N}$\footnote{Note that 
 the number of particles does not have to be $N$ but we make this assumption
 for notational convenience.}
 by employing a local regression within the clusters associated 
 with the $S_k^\ast$'s defined above. More precisely, we write
\begin{align}\label{eq:weighted-post}
     T_t(x,y) 
     &\approx \widehat T_t \rebuttal{(x,y)}:=\sum_{k=1}^K w_t^k\, T_k^\ast(x,y),\\
     w_t^k 
     &:= \exp\Bigl(-\lambda\, \rho(S_k^\ast,S_t)\Bigr)/\sum_{k'=1}^K \exp\Bigl(-\lambda\, \rho(S_{k'}^\ast,S_t)\Bigr),\nonumber
\end{align}
where $S_t$ is defined as  $\Bigl((X^i_{t|t-1},Y^i_{t|t-1})_{i=1}^{N}, 
T_t\Bigr)$, $T^\ast_k$ is the OTF map associated with the medoid 
$S^\ast_k$, $\lambda \geq 0$ is a hyperparameter, and $\rho$ is a distance function similar to $d$. We introduce $\rho$ here to distinguish between offline and online distance functions and analogously we use similar subscript notations (i.e. $\rho_{W_2},\rho_{MMD}$) to refer to the choice of the distance function. Since $T_t$ is unknown, the distance $\rho(S_k^\ast,S_t)$ is computed solely using $X$-samples. 

\begin{remark} 
There is a trade-off between accuracy and online 
computational cost within the choice of $K$, with larger values of $K$ slowing down 
computations (since at each inference step one must evaluate distances 
between $S_t$ and the $K$-medoids $S^\ast_k$) while improving 
accuracy since more maps will be used to construct $T_t$.
\end{remark}
 



\begin{remark}
     One can compute $d_{W_2},d_{MMD}$ using the samples from the joint distribution by concatenating $(X^i,Y^i)$. However, in our numerical experiments, this approach did not yield a significant improvement in performance.
\end{remark}



\subsection{Theoretical Intuition for A-OTF}
\rebuttal{
The theoretical justification for our proposed interpolation~\eqref{eq:weighted-post} comes from theory of non-parametric estimation~\cite{tibshirani2023nonparametric}. In particular, the case $\lambda=\infty$ in~\eqref{eq:weighted-post} corresponds to the nearest-neighbor estimation, while finite values of $\lambda$ corresponds to a  kernel smoothing estimator. The theory of non-parametric estimation implies consistency  under regularity assumption for the operator that is being estimated with respect to the metric $\rho$. 
}
Recent works in the theory of Bayesian inference \cite{sprungk2020local, garbuno2023bayesian} 
have shown that posterior distributions 
vary continuously as  a function of their priors 
when these perturbations are measured with respect 
to OT-type metrics. Similarly, recent works on  conditional OT \cite{hosseini2025conditional,chemseddine2024conditional}
further show that the maps $T$, that characterize 
the Bayesian updates in our setting, vary continuously 
with respect to the priors and the observations. 
These results imply that, under sufficient assumptions, 
the mapping $( X^i_{t|t-1}, Y^i_{t|t-1})_{i=1}^N  \mapsto T_t$, 
viewed as a mapping from empirical distributions to 
transport maps, is continuous and hence amenable to 
numerical approximation. Thus, we conjecture that
if the A-OTF training data is sufficiently ``space 
filling" or ``dense" then the resulting A-OTF 
approximation should converge to the true solution 
of \eqref{eq:new_loss}. A proof of this conjecture is 
left as a future research direction and not investigated here.


    


\section{Numerical experiments}\label{sec:numerics}
\subsection{The numerical algorithm}\label{sec:numerical_alg}

We summarize the algorithmic 
implementation of A-OTF. As mentioned earlier, the online stage concerns the numerical 
implementation of the $K$-medoids algorithm with the 
main computational cost coming from solving the 
individual OTF problems to construct the A-OTF 
training data. For this we used the OTF implementation in 
\cite{al-jarrah2024nonlinear}. The computation of the distance matrix
$D$ was also done using off-the-shelf OT solvers for 
the case of $d_{W_2}$, while $d_{MMD}$ 
was implemented with the radial basis function kernel (RBF), normalizing the data, and tuning the length scale to a fixed value among all experiments. A summary of the A-OTF procedure appears in Algorithm~\ref{alg:A-OTF}.






\begin{algorithm}[t]
\caption{A-OTF algorithm} 
\begin{algorithmic}
\STATE \hspace{-10pt}\textbf{Offline Stage:}
\STATE \textbf{Input:} Training data $\{\widetilde S_m\}_{m=1}^M$, number of clusters $K$, the choice of distance function $d$. 
\STATE \textbf{Distance matrix:} Compute $D$ according to~\eqref{eq:Dmn}.

\STATE \textbf{Clustering:} Use $K$-Medoids algorithm with the distance matrix $D$ to select $K$ representatives $\{S_k^\ast\}_{k=1}^K$.

\STATE \textbf{Output:} $\{S_k^\ast\}_{k=1}^K$.\\[7pt]

\STATE \hspace{-10pt}\textbf{Online Stage:}
\STATE \textbf{Input:} $\{S_k^\ast\}_{k=1}^K$, $\{ X^i_{0|0}\}_{i=1}^N \sim \pi_0$, distance function $\rho$.
\FOR{$t=1$ to $t_f$}
\STATE \textbf{Propagation:} Update $X_{t|t-1}^i$ according to~\eqref{eq:X_particle_law} and $Y_{t|t-1}^i\sim h_t(\cdot|X_{t|t-1}^i) \quad \forall i=1,\dots,N.$ 
\STATE \textbf{Compute $\widehat{T}_t$:} Compute $\{w_t^k\}_{k=1}^K$ and $\widehat{T}_t$ as~\eqref{eq:weighted-post}.
\STATE \textbf{Conditioning:} $X_{t|t}^i = \widehat{T}_t(X_{t|t-1}^i,Y_t)\quad \forall i=1,\dots,N.$
\ENDFOR

\STATE \textbf{Output:} $\{X_{t|t}^i\}_{i=1}^N \quad \textit{for} ~ t=1,\dots,t_f.$
\end{algorithmic}
\label{alg:A-OTF}
\end{algorithm}

\subsection{Experiment Setup}

We compared A-OTF against three other algorithms for filtering: the Ensemble Kalman filter (EnKF)~\cite{evensen2006}, the 
SIR particle filter~\cite{doucet09}, and OTF~\cite{al-jarrah2024nonlinear}. The details of all three algorithms appear in~\cite{al-jarrah2024nonlinear}, and the numerical code used to produce the results is available online\footnote{\url{https://github.com/Mohd9485/A-OTF}}. 


\subsection{Lorenz 63}\label{sec:L63}
\rebuttal{
For our first experiment, we consider the  Lorenz 63~\cite{lorenz1963deterministic} with the observation model and the prior distribution:
\begin{equation}\label{eq:L63}
Y_t = X_t(3)
+ \sigma_{obs}W_t, \quad X_{0} \sim \mathcal{N}(\mu_0,\sigma_0^2I_3),
\end{equation}
where \([X(1),X(2),X(3)]^\top\) represent the hidden states of the system,
$\mu_0 = [0,0,0]^\top$, and $\sigma_{0}^2=10$. 
We used a time-discretization of the dynamics, with $\Delta t = 0.01$, to bring the model into our setup~\eqref{eq:model-dyn}. 
The noise $W_t$ is a standard Gaussian random variable and $\sigma_{obs}^2=10$. Observing the third state $X(3)$ results in bimodality of the posterior in the first two states due to the symmetry in the Lorentz 63 dynamic model.

\begin{figure}[t]

         \centering
         \includegraphics[width=0.485\hsize,trim={15 0 50 60},clip]{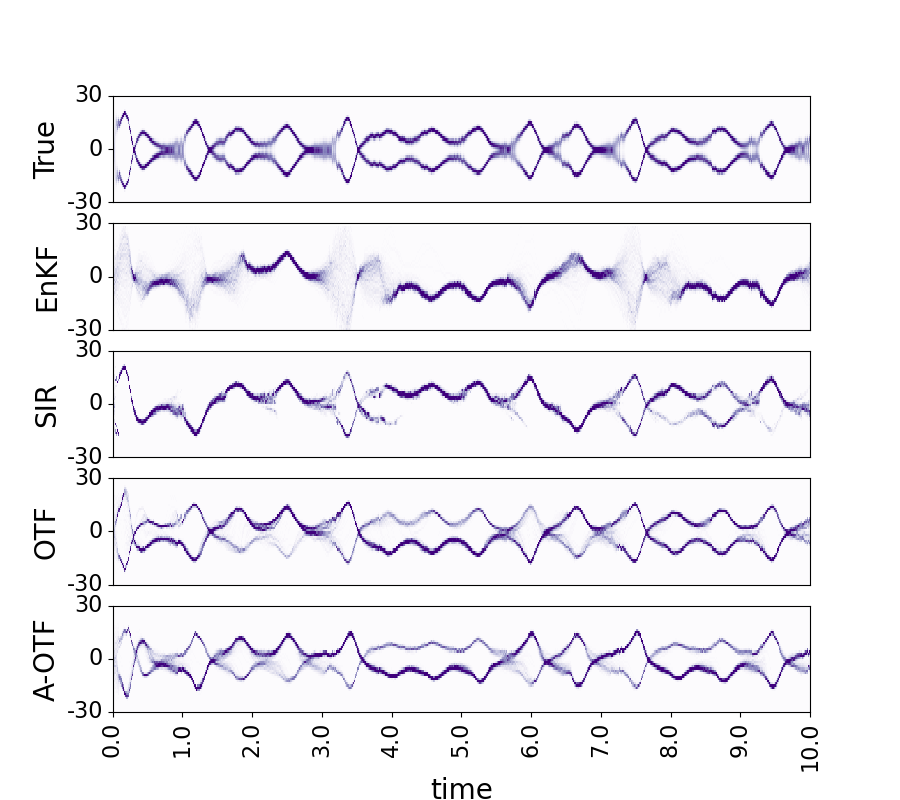}
         \hspace{0.001cm}\includegraphics[width=0.485\hsize,trim={15 0 50 60},clip]{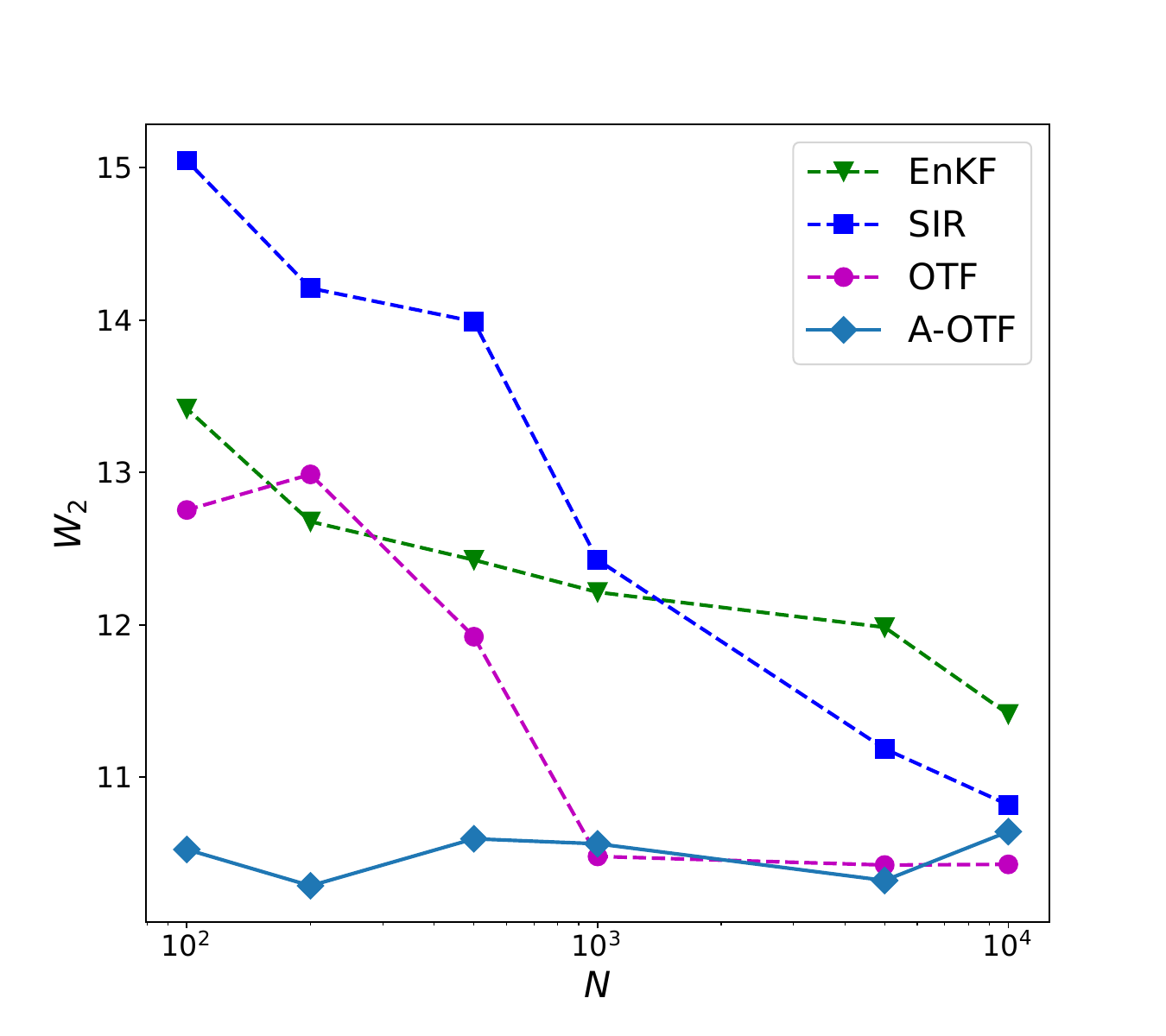}

     \caption{
     Numerical results for the Lorenz 63 example, in Section~\ref{sec:L63}. The left column shows the particle trajectory distributions of the first unobserved components of the particles along with the true distribution, where $A-OTF$ uses $d_{W_2},\rho_{W_2}$, $K=20$, and $\lambda=1$. The right column presents the empirical $W_2$ distances between each method and the true distribution as a function of the number of particles $N$ for a fixed $\mu_0 = 4$ and $\sigma_0=5$ averaged over five independent simulations and each simulation consists of 500 time steps.  
     }
    \label{fig:L63}
\end{figure}
We implemented the offline stage of A-OTF with $N=1000$ and $M=1000$ to create the dataset of  pre-trained maps and used $d_{W_2}$ to perform $K$-medoids with $K=20$. We implemented the online stage of A-OTF with $\rho_{W_2}$ and $N=250$, which is the same for all other methods. The result is depicted in  
Figure~\ref{fig:L63}. The left panel shows the distribution of the first component of the particle trajectories, along with the true posterior distribution computed by simulating the SIR with $10^5$ particles. The bimodality of the true posterior is replicated in OTF and A-OTF, while EnKF and SIR fail due to small number of particles. A quantitative evaluation of the error, as the number of particles $N$ increases, is depicted in the right panel. The error is evaluated by computing the empirical Wasserstein distance between the  particle distributions of each algorithm, and the ground-truth SIR with $10^5$ particles. The result shows that the performance of the SIR  becomes as good as A-OTF and OTF when $N$ becomes large. Therefore, due to the additional computational cost of A-OTF, it seems that, given the same computational budget, the performance of the SIR is better than A-OTF.  
We study this question in the next section and show that this observation is not valid when the dimension of the problem becomes larger. 
}

\rebuttal{
\subsection{Linear dynamics with quadratic observation}\label{sec:lin-example}
Next, we consider
the linear dynamics and quadratic observation model 
\begin{equation}\label{eq:model-example}
\begin{split}
    X_{t} &= A
    X_{t-1} + \sigma V_t,~
    Y_t = X_t^2 + \sigma W_t
\end{split}
\end{equation}
where $A = \mathrm{diag}(F, F, \dots, F) \in \mathbb{R}^{2n\times 2n}$, 
\[
F = \begin{bmatrix}
        \alpha & \sqrt{1-\alpha^2}
        \\
        -\sqrt{1-\alpha^2} & \alpha
    \end{bmatrix}
\] 
$X_t\in \mathbb R^{2n}$, $Y_t \in \mathbb R^{2n}$, $\{V_t\}_{t=1}^\infty$ and $\{W_t\}_{t=1}^\infty$ are i.i.d sequences of $2n$-dimensional standard Gaussian random variables, $\alpha=0.9$, and $\sigma^2=0.01$. The quadratic observation model leads to a bi-modal posterior. We used this model as a validation test, with $n=1$ (i.e. two-dimensional) as a low dimensional setting and $n=4$ (i.e. eight-dimensional)  as a high dimensional setting. 

The numerical results for $n=1$ are presented in Figure~\ref{fig:XX_n_1}, which follows the same structure as Figure~\ref{fig:L63}. We implemented the offline stage of A-OTF with $N=10^4$ and $M=500$ to create the dataset of  pre-trained maps and used $d_{W_2}$ to perform $K$-medoids with $K=10$. We implemented the online stage of A-OTF with $\rho_{W_2}$ and $N=10^4$, which is the same for all other methods. The ground-truth is provided by SIR with $N=10^5$. The results confirm our observations from Figure~\ref{fig:L63}.


The numerical results for $n=4$ are presented in Figure~\ref{fig:XX_n_4} and Figure~\ref{fig:XX_n_4_change_IC_N}. 
The implementation of A-OTF and other algorithms are the same as the ones in $n=1$. For ground-truth, we utilized the independence structure of the filtering problem and simulated $n=4$ independent SIR algorithms with $N=10^5$ for each independent component pair. Notably, the left panel in Figure~\ref{fig:XX_n_4} shows that SIR fails to capture the bimodality in this higher dimensional setting. 

To study the effect of the distance metrics, we implemented different variations of algorithm~\ref{alg:A-OTF} with different distance functions $d_T,d_{W_2},d_{MMD}$ to compute the distance matrix $D$. Then we used different distance functions for $\rho_{W_2},\rho_{MMD}$ in the online stage where we refer to setting $\lambda=1$ as \emph{Weighted} and $\lambda=\infty$ as \emph{Nearest}. The right three columns of Figure~\ref{fig:XX_n_4} show the $W_2$ distance, between the output of each variation and the true distribution, as a function of the number of selected maps $K$. 
The results indicate all methods have similar trends with a slight increase in performance and robustness for the weighted methods across different distance functions $d$.  



We provide numerical evidence of the robustness of the A-OTF by fixing the pre-trained maps used in the offline stage while varying the mean $\mu_0$ and variance $\sigma_0$ of the prior distribution in~\eqref{eq:model-example} that is used to create the observational data for the online stage. The left two panels of Figure~\ref{fig:XX_n_4_change_IC_N} illustrate the empirical $W_2$ distance between the true posterior and the output of each method, as a function of the parameters $\mu_0$ and $\sigma_0$ where all methods use $N=10^4$ particles. The results indicate that A-OTF is robust under the change in $\mu_0$ and $\sigma_0$ and provides a significant improvement over SIR and EnKF methods and a slight improvement over the regular OTF method. 

In order to study the effect of the number of particles, we performed the same experiment with varying $N$ while fixing $\mu_0=0$ and $\sigma_0=1$. It is observed that OTF and A-OTF methods perform better than other methods. In this example, SIR requires a large number of samples $(\gg 10^6)$ since the system is $8$ dimensional and the likelihood is highly degenerate due to the small noise level $\sigma^2$.  
Finally, the fourth panel shows the same error metric as a function of computational time (as $N$ increases). The result illustrates a clear evidence for the  computational efficiency of A-OTF in comparison to the other methods, in this higher dimensional setting.  It is also noted that additional computational savings are possible by following more efficient computations of distance functions (e.g. by sub-sampling and using Sinkhorn-type distances instead of $W_2$). 

\begin{figure}[t]

         \centering
         \includegraphics[width=0.485\hsize,trim={15 0 50 60},clip]{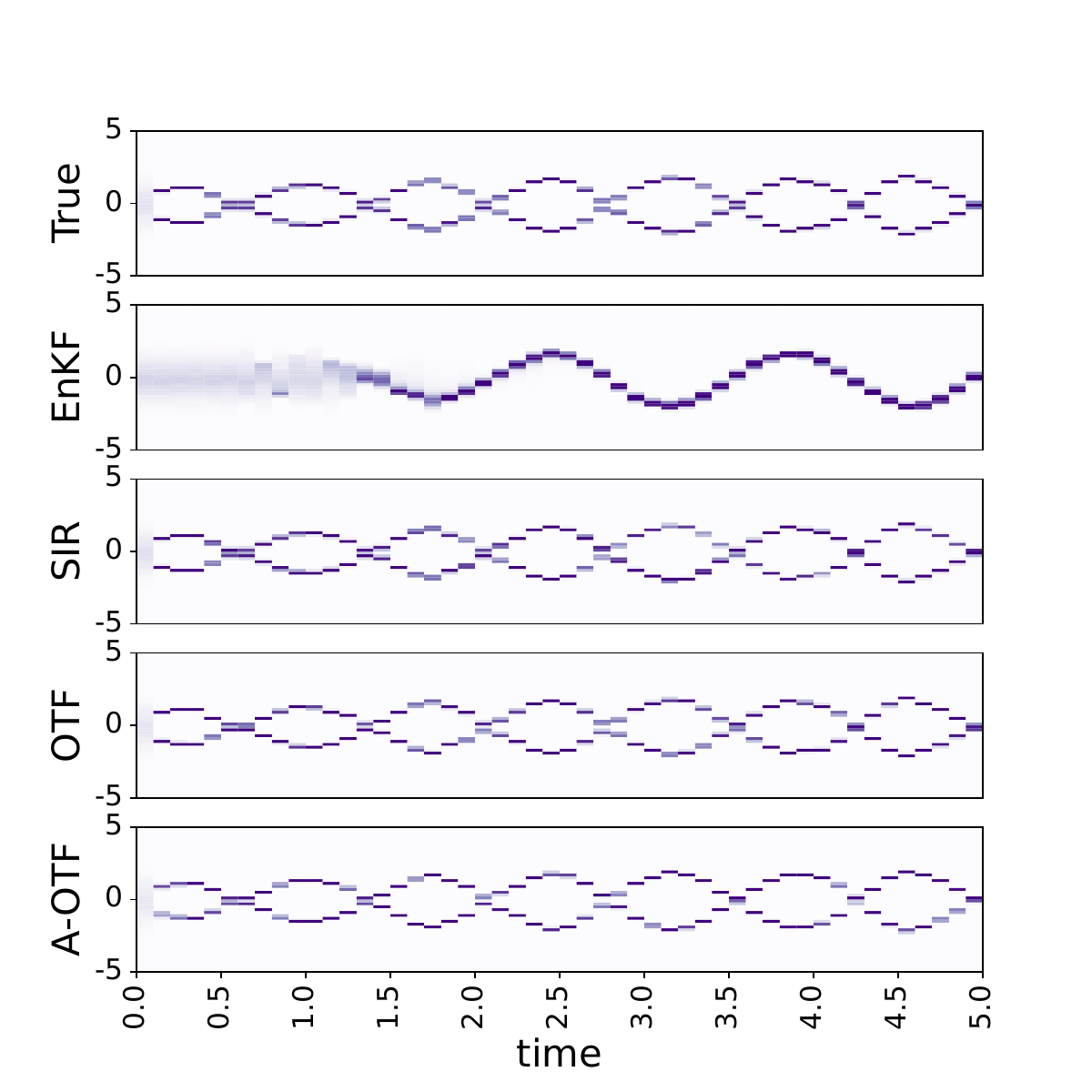}
         \hspace{0.001cm}
         \includegraphics[width=0.485\hsize,trim={15 0 50 60},clip]{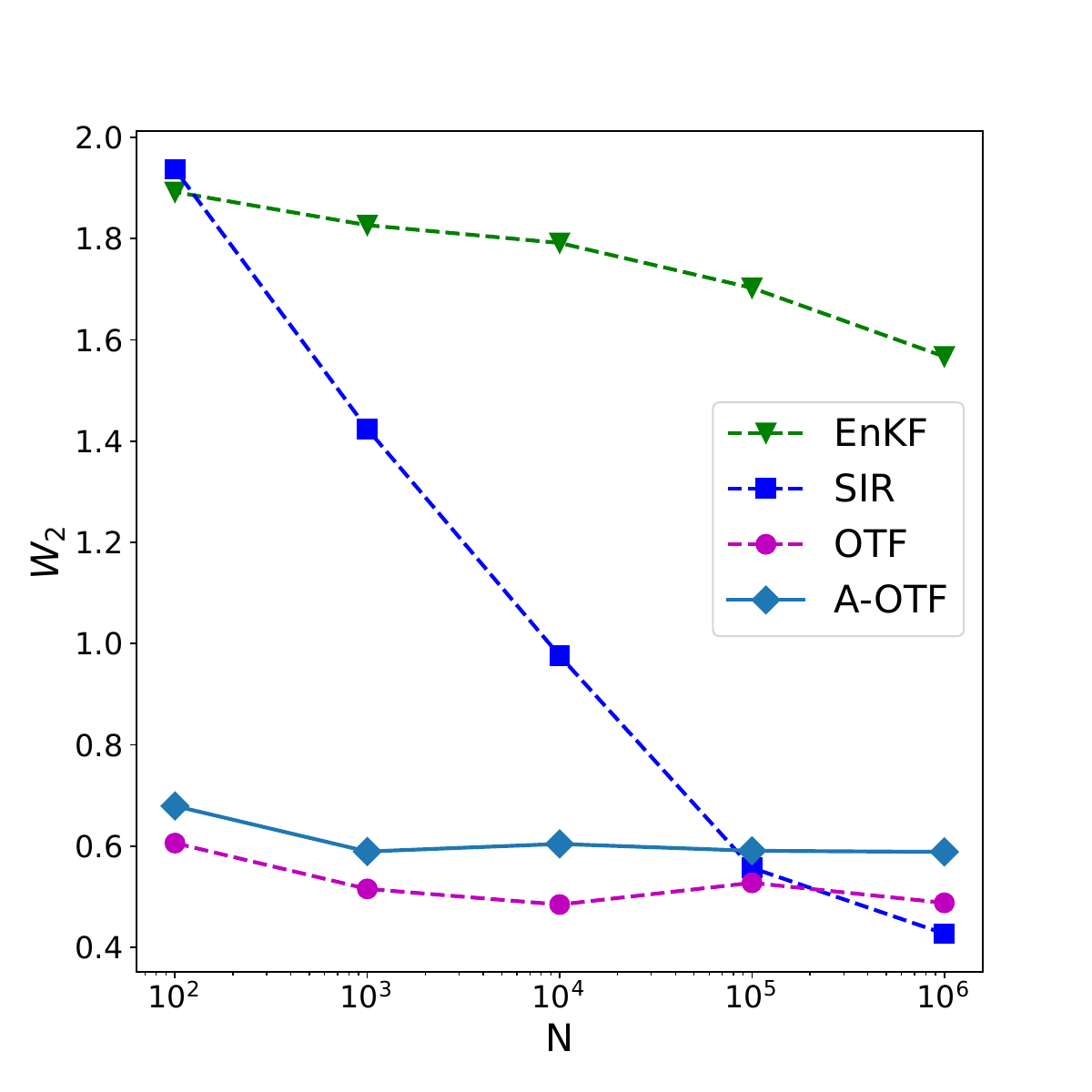}
     \caption{
        Numerical results for the quadratic observation  example in Section~\ref{sec:lin-example}. The left panel shows the particles distribution for different methods compared to the true distribution as a function of time. The right panel shows the empirical $W_2$ distance between the true distribution and the output distribution of each method, averaged over five independent simulations as a function of the number of particles $N$.
     }
    \label{fig:XX_n_1}
\end{figure}

\begin{figure*}[ht]

         \centering
         \includegraphics[width=0.24\hsize,trim={0 0 53 60},clip]{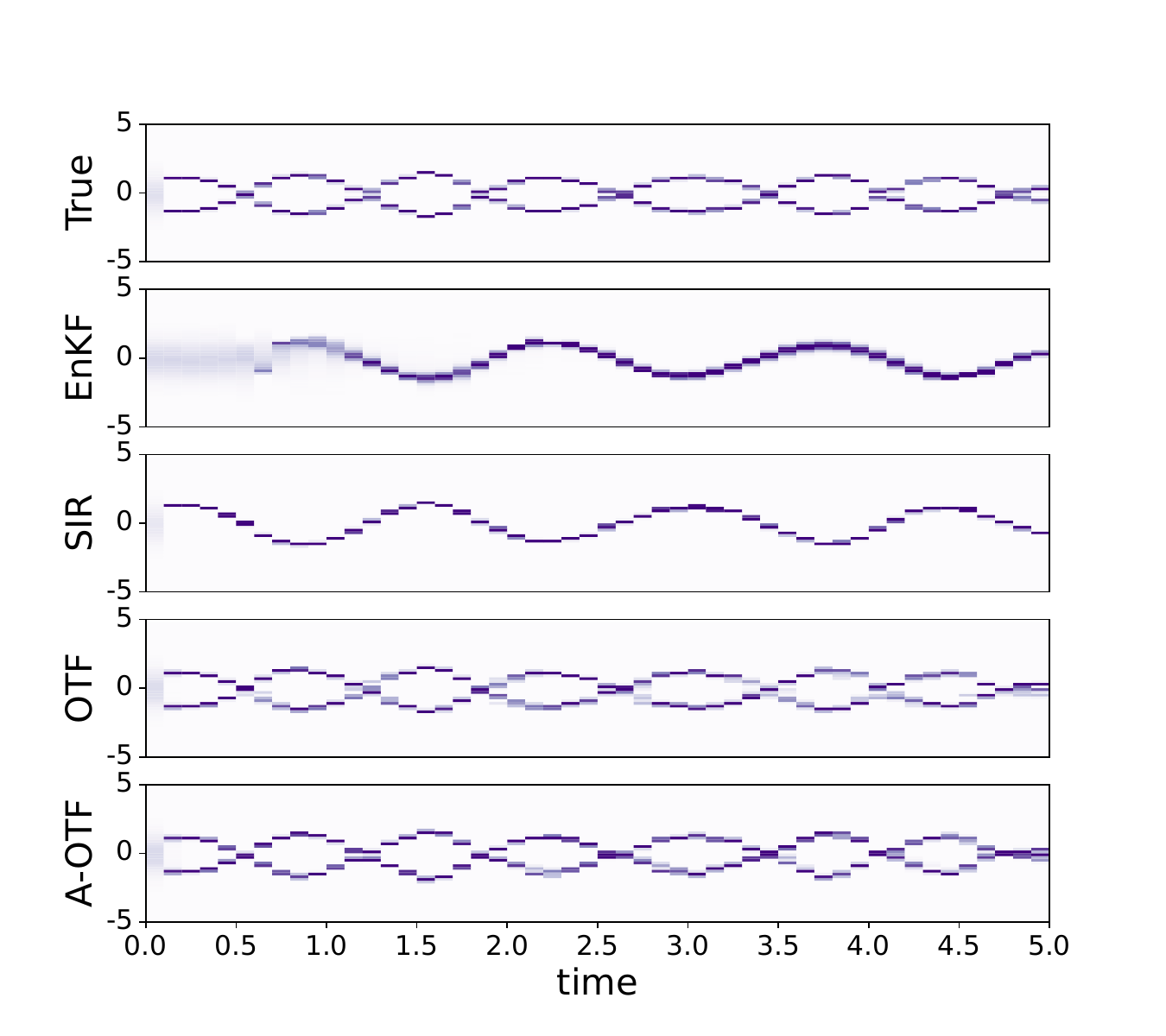}
         \hspace{0.001cm}
         \includegraphics[width=0.24\hsize,trim={0 0 53 60},clip]{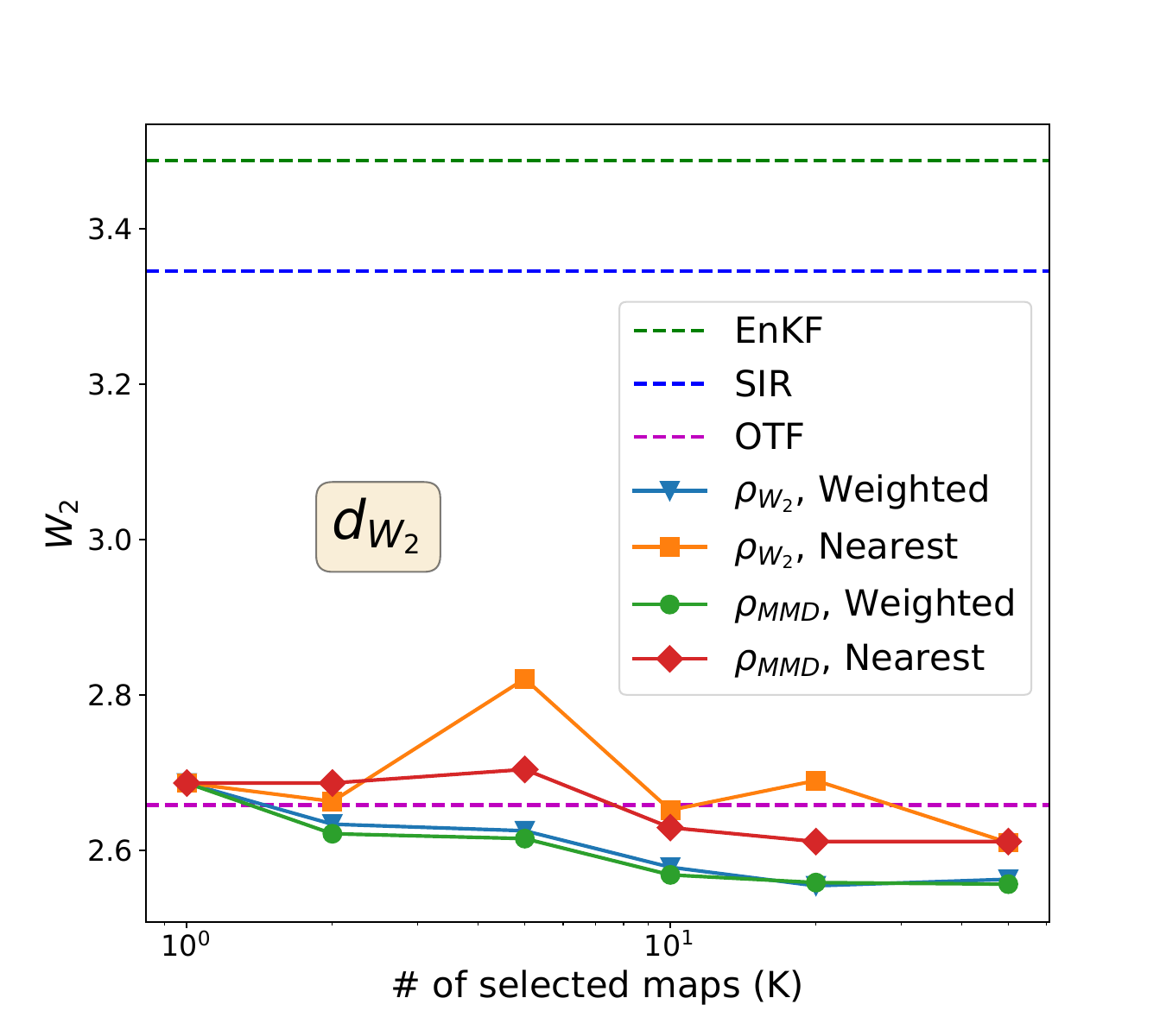}
         \hspace{0.001cm}
         \includegraphics[width=0.24\hsize,trim={0 0 53 60},clip]{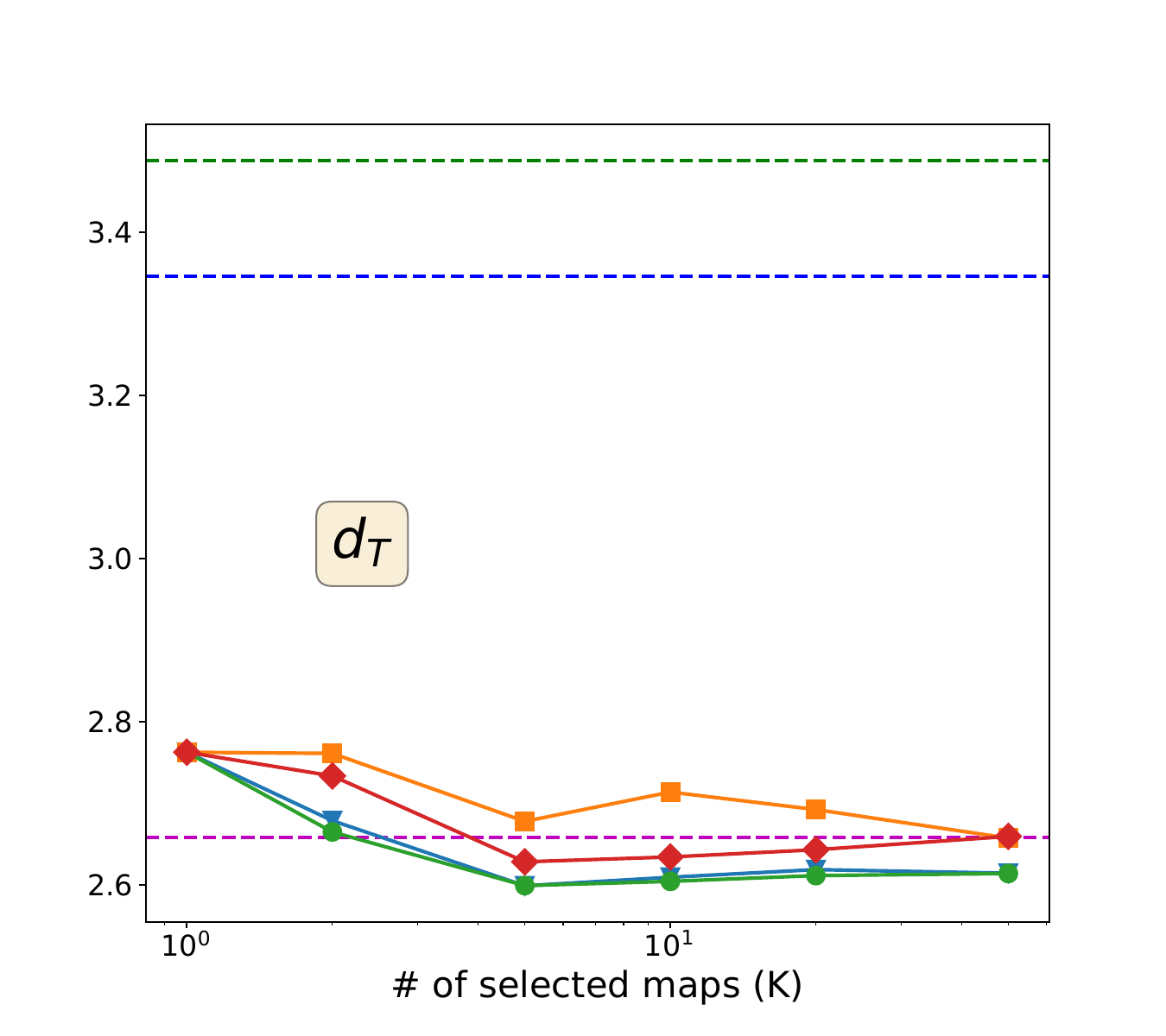}
         \hspace{0.001cm}
         \includegraphics[width=0.24\hsize,trim={0 0 53 60},clip]{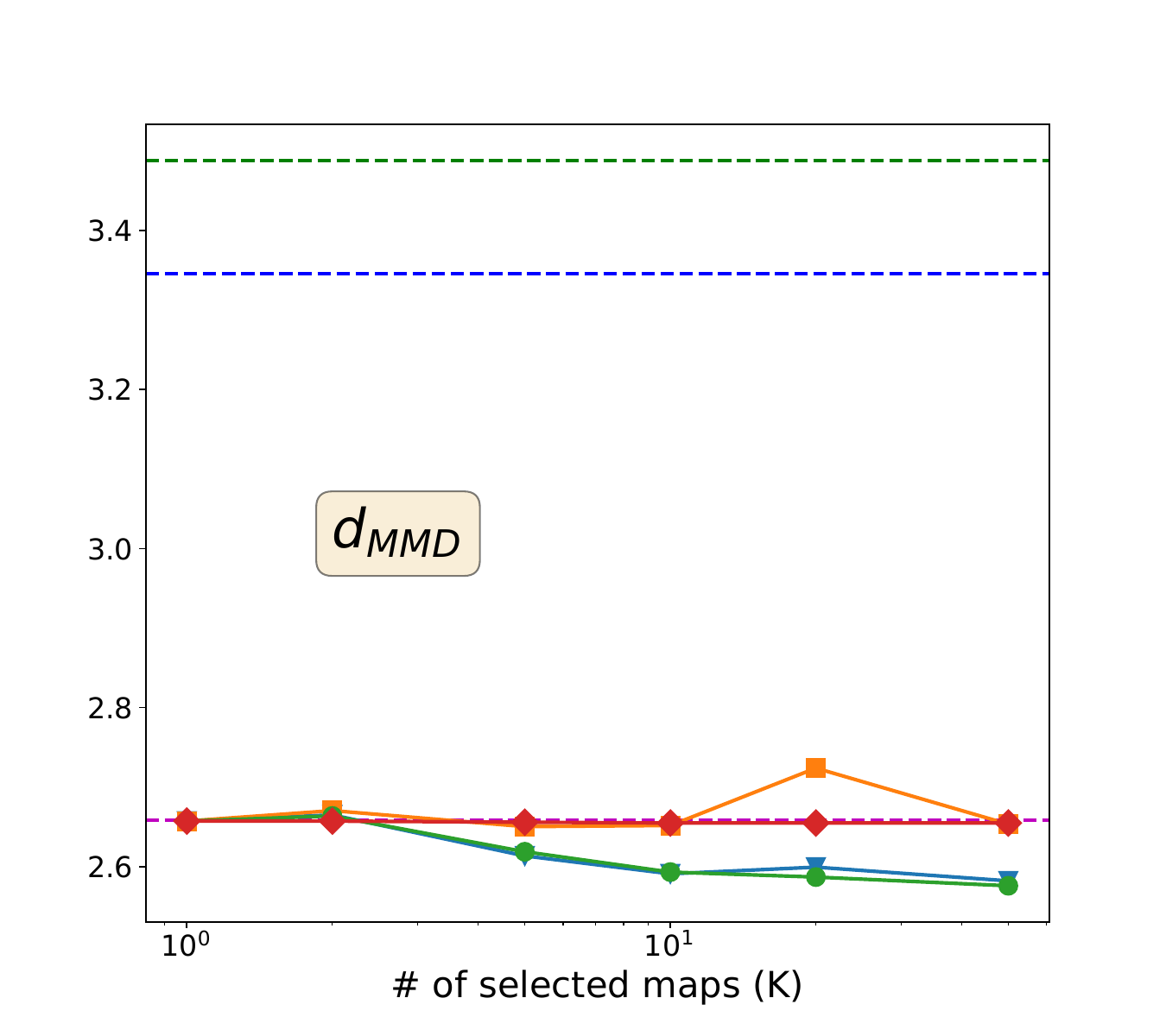}
     
     \caption{
        Numerical results for the quadratic observation example in Section~\ref{sec:lin-example}. The left column shows the particle trajectory distributions of the first unobserved components of the particles along with the true distribution.  The right three columns show the empirical \(W_2\) distance between the true distribution and the output distribution of each method, averaged over five independent simulations and $50$ time steps, as a function of the number \(K\) of selected maps according to \(d_{W_2\!}\) (left panel), \(d_T\) (middle panel), and \(d_{MMD\!}\) (right panel), respectively. 
     }
    \label{fig:XX_n_4}
\end{figure*}

\begin{figure*}[ht]

         \centering
         \includegraphics[width=0.24\hsize,trim={0 0 53 60},clip]{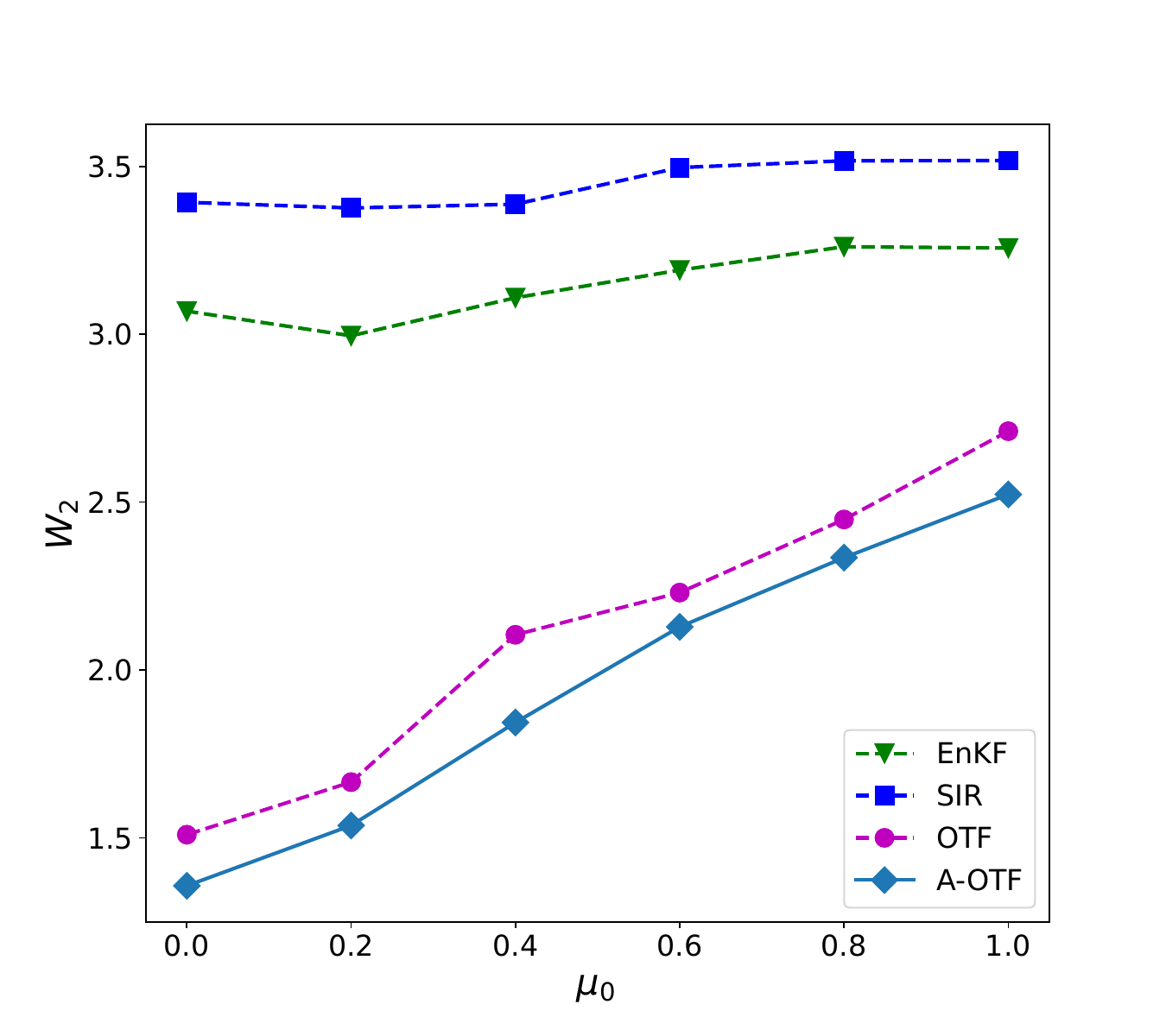}
         \hspace{0.001cm}
         \includegraphics[width=0.24\hsize,trim={0 0 53 60},clip]{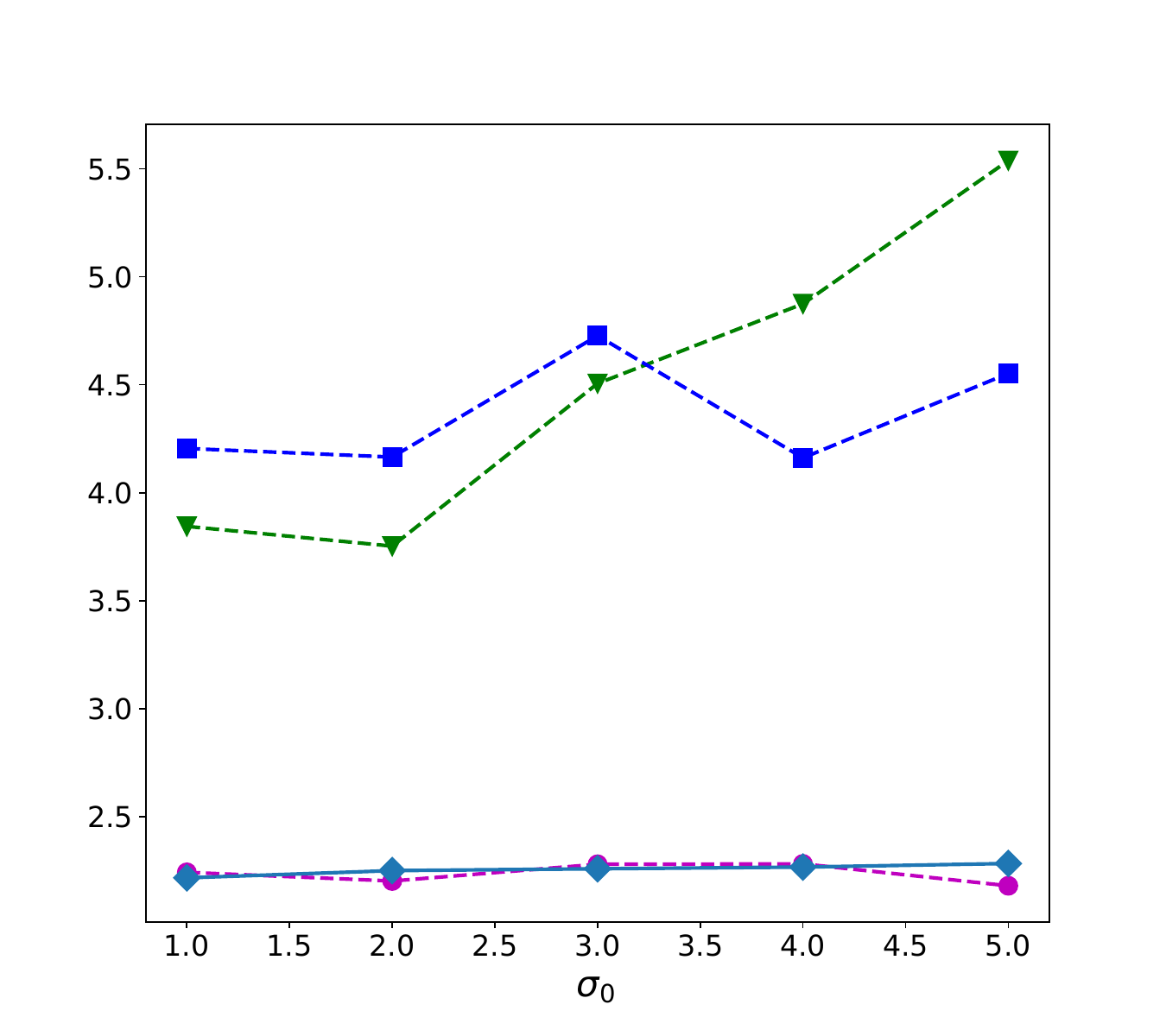}
         \hspace{0.001cm}
         \includegraphics[width=0.24\hsize,trim={0 0 53 60},clip]{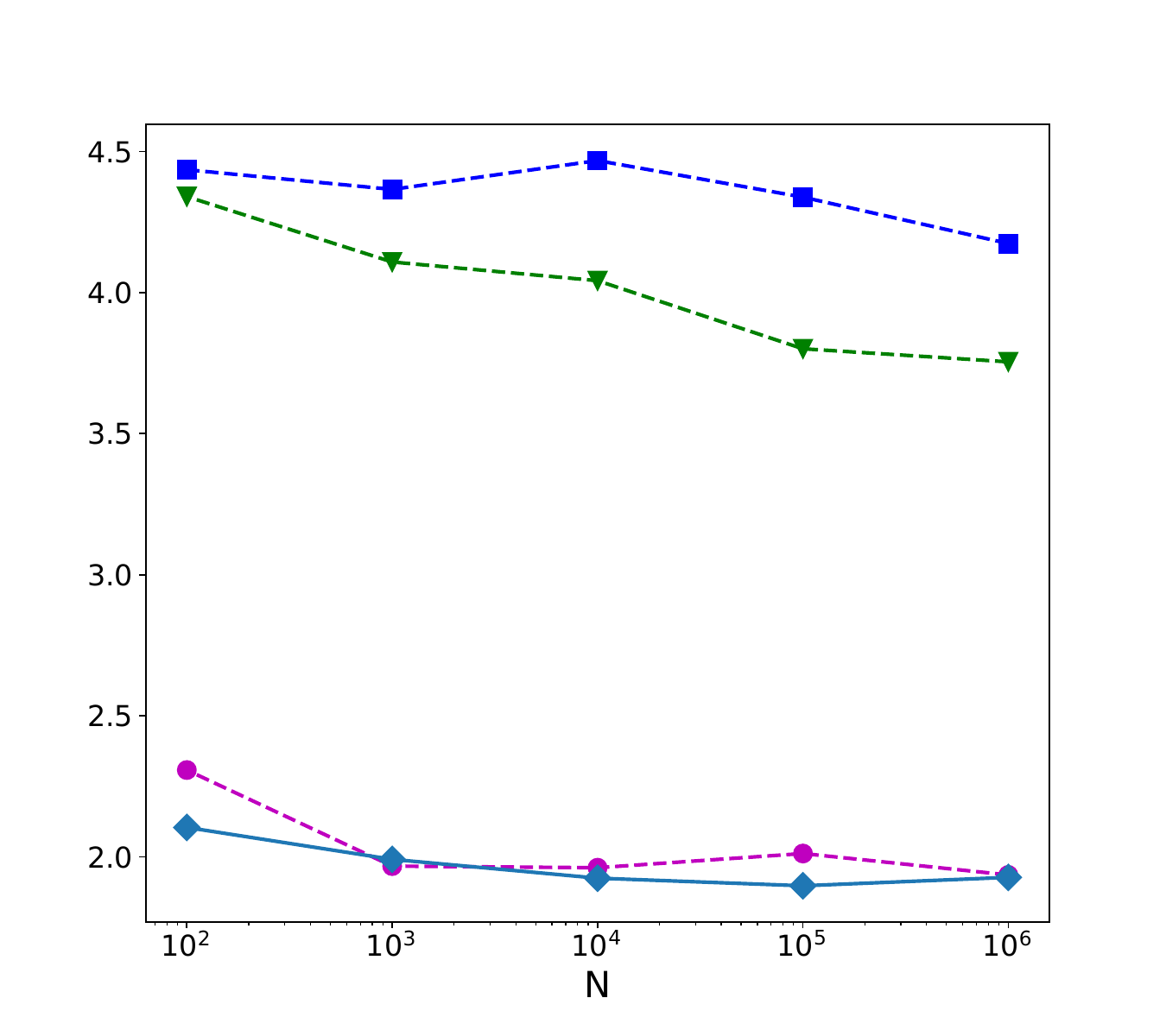}
         \hspace{0.001cm}
         \includegraphics[width=0.24\hsize,trim={0 0 53 60},clip]{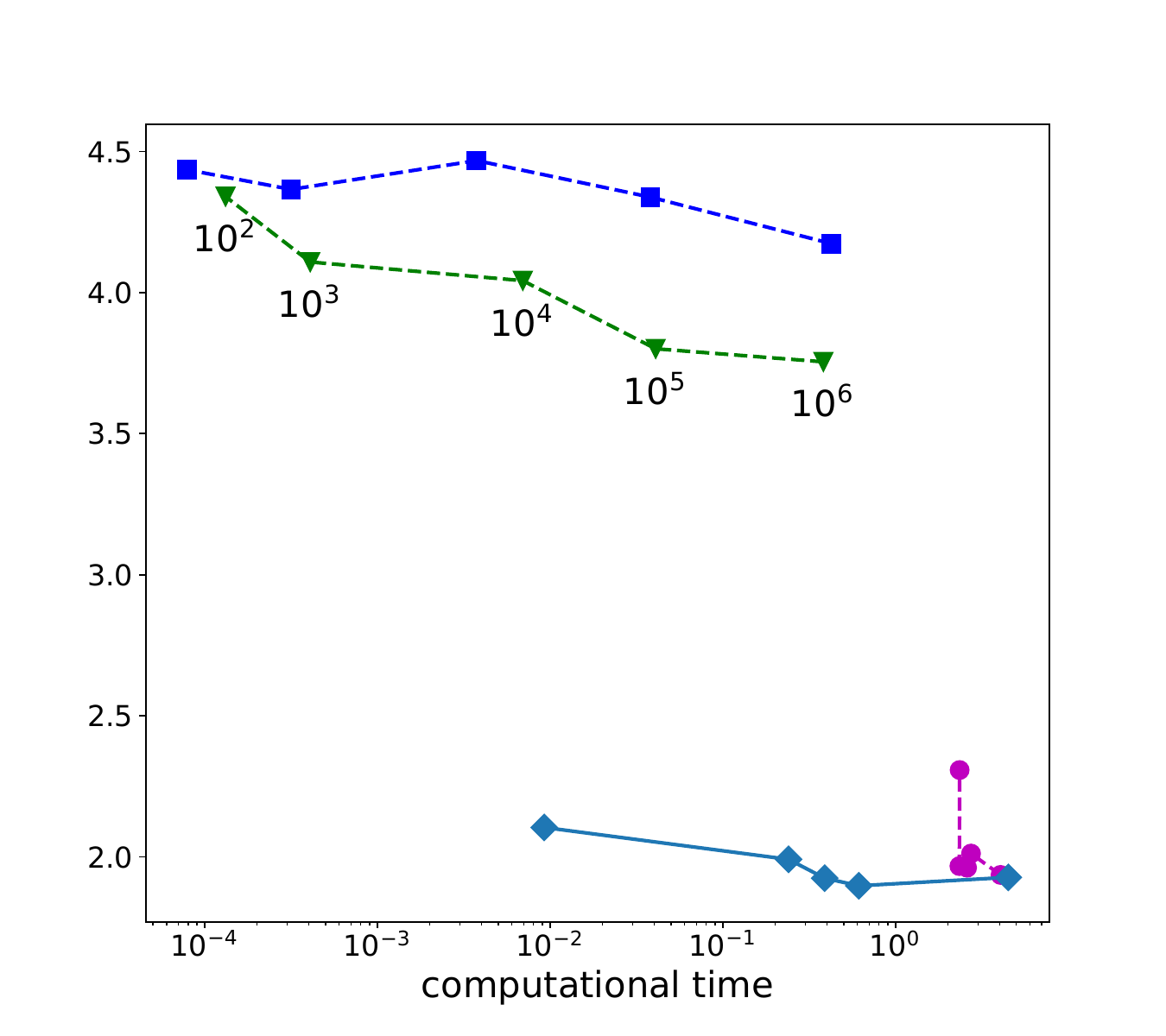}
     
     \caption{
        Numerical results for the quadratic observation example in Section~\ref{sec:lin-example}. The first two panels present the empirical \(W_2\) distances between each method and the true distribution as a function of the particle mean \(\mu_0\) and variance \(\sigma_0\), while the amortized pre-trained maps are fixed. The third and fourth panels present the same error metric as as a function of the number of particles $N$ and the computational time, respectively.   All results are averaged over five independent simulations. 
     }
    \label{fig:XX_n_4_change_IC_N}
\end{figure*}
}

\section{Discussion}\label{sec:discussion}
In this work, we introduced A-OTF, a novel filtering approach that leverages pre-trained OTF maps and capitalizes on their inherent similarities during an offline phase to mitigate the computational cost associated with online inference. The framework eliminates the need for online training, thereby significantly reducing computational overhead and avoiding issues associated with hyperparameter tuning. Consequently, A-OTF represents a significant step towards developing computationally efficient and practically viable variants of OTF algorithms. Future work will focus on extending the mathematical justification and error analysis of this approach, while also exploring its adaptability across a broader range of applications.

    
    \bibliographystyle{plain}
    \bibliography{references}

\end{document}